\documentclass{article}

\usepackage{graphicx}
\usepackage{amsmath}
\usepackage{amsfonts}
\usepackage{amssymb}

\usepackage{bm}

\usepackage{graphicx}
\usepackage{amsmath}
\usepackage{amsfonts}
\usepackage{amssymb}

\newcommand*{\de}{\mathop{}\!\mathrm{d}}

\def\ba{\begin{aligned}}

\def\be{\begin{equation}}

\def\cl {\nonumber \\}

\def\ea{\end{aligned}}

\def\ee{\end{equation}}

\def\el {\nonumber }

\def\g{{\bm g}}

\def\n{{\bm n}}

\begin{document}

\title{Data-Driven Enhanced Model Reduction for Bifurcating Models in Computational Fluid Dynamics}

\author{Martin W. Hess\footnote{Martin W. Hess, SISSA mathLab, International School for Advanced Studies, via Bonomea 265, I-34136 Trieste, Italy, \texttt{ mhess@sissa.it}}, 
Annalisa Quaini\footnote{Annalisa Quaini, Department of Mathematics, University of Houston, Houston, Texas, 77204, USA, \texttt{quaini@math.uh.edu}} , 
and Gianluigi Rozza\footnote{Gianluigi Rozza, SISSA mathLab, International School for Advanced Studies, via Bonomea 265, I-34136 Trieste, Italy, \texttt{grozza@sissa.it}}}


\maketitle

\abstract{We investigate various data-driven methods to enhance pro\-jection-based model reduction techniques with the aim of capturing bifurcating solutions. To show the effectiveness of the data-driven enhancements, we focus on the incompressible Navier-Stokes equations and different types of bifurcations. 
To recover solutions past a Hopf bifurcation, we propose 
an approach that combines proper orthogonal decomposition
with Hankel dynamic mode decomposition. To approximate solutions close to a pitchfork bifurcation, we combine localized reduced models with artificial neural networks.
Several numerical examples are shown to demonstrate
the feasibility of the presented approaches.
}

\section{Introduction}

Model Order Reduction is a vital tool for parametric studies of systems governed by partial differential equations as it drastically cuts the computational time required by standard approximations, such as the Spectral Element method or the Finite Element method. Data-Driven methods can be used to enhance existing intrusive reduced order models (ROMs) with the aim of
applying them to problems previously deemed overly complicated and achieving greater accuracy. 

This work focuses on bifurcating models arising in incompressible fluid dynamics. Specifically, we investigate a jet flow that exhibits a supercritical pitchfork bifurcation under a change in the Reynolds number and jet inlet width and a cavity flow that shows multiple Hopf bifurcations under a change of Grashof number. Bifurcating models pose particular challenges for model order reduction. When a pitchfork bifurcation is present in the parameter domain of interest, multiple solutions exist
in certain regions of the domain while in other regions only one solution exists. Transitions from steady solutions to time-periodic and finally fully time-dependent solutions without any discernible period occur as a parameter is varied in a domain containing Hopf bifurcations. 
Various intrusive and also data-driven ROMs for bifurcating problems 
have also been investigated in \cite{pichi2020optcntrl, PITTON2017534, pichi2021artificial,https://doi.org/10.1002/fld.5118}. 

In this paper, we discuss several data-driven enhancements to existing ROMs. We focus on two variants of the dynamic mode decomposition (DMD) algorithm, namely the Hankel-DMD and a stabilized DMD, to improve the ROM for the cavity
flow problem introduced in \cite{HessQuainiRozza2022_ETNA}. 
Related work using the DMD can be found in \cite{GADALLA2021104819}.
For the jet flow problem, we propose to use
artificial neural networks (ANNs) to assign a local ROM
\cite{Hess2019CMAME} to a given region of the parameter domain. This latter approach ultimately enables the recovery of multiple solutions at a single parameter point. All of the above-mentioned 
enhancements have been implemented in the open source software ITHACA-SEM \cite{HessLarioMengaldoRozza_ICOSAHOM_2022} and are available to the ROM community.

The remainder of the paper is structured as follows. 
Section 2 introduces the abstract model setting and specifies it for the jet and cavity flow problems. Section 3 and 4 discuss the model reduction approach for the jet and cavity problems, respectively. In Section 5, we draw conclusions and provide further perspectives.

\section{Application to the incompressible Navier-Stokes equations}

The motion of an incompressible, viscous, and Newtonian fluid in a spatial domain $\Omega \subset\mathbb{R}^d$, $d = 2$ or $3$, over a time interval of interest $(0, T]$ is governed by the Navier-Stokes equations:
\begin{eqnarray}\label{NS-1}
\frac{\partial {\bm u}}{\partial t}+({\bm u}\cdot  \nabla {\bm u})-\nu\Delta {\bm u}+\nabla p={\bm f} &\qquad\mbox{in } \Omega\times(0,T],\\ \label{NS-2}
\nabla \cdot {\bm u}=0&\qquad\mbox{in } \Omega\times (0,T].
\end{eqnarray}
Here, ${\bf u}$ and $p$ are the unknown velocity and pressure fields, $\nu>0$ is the kinematic viscosity of the fluid, and ${\bm f}$ denotes the body forces. 
If the system has evolved towards a steady state, then $\frac{\partial {\bm u}}{\partial t}$
in \eqref{NS-1} can be neglected.

Initial and boundary conditions complete problem \eqref{NS-1}--\eqref{NS-2}
\begin{eqnarray}
{\bm u}={\bm u}_0&\qquad&\mbox{in } \Omega \times \{0\} \label{IC}\\
{\bm u}={\bm u}_D&\qquad&\mbox{on } \partial\Omega_D\times (0,T] \label{BC-D} \\
-p \n + \nu \frac{\partial {\bm u}}{\partial \n}=\g &\qquad&\mbox{on } \partial\Omega_N\times, (0,T], \label{BC-N}
\end{eqnarray}
where $\partial\Omega_D \cap \partial\Omega_N = \emptyset$ and $\overline{\partial\Omega_D} \cup \overline{\partial\Omega_N} =  \overline{\partial\Omega}$. Here, ${\bm u}_0$, ${\bm u}_D$, and $\g$ are given and $\n$ denotes the unit normal vector on the boundary $\partial\Omega_N$ directed outwards. The boundary $\partial \Omega_D$ is called the Dirichlet boundary and $\partial \Omega_N$ the Neumann boundary.

Let $L^2(\Omega)$ denote the space of square integrable functions in $\Omega$ and $H^1(\Omega)$ the space of functions belonging to $L^2(\Omega)$ with weak first derivatives in $L^2(\Omega)$. 
We define sets 
\begin{eqnarray}
{\bm V} &:=& \left\{ {\bm v} \in [H^1(\Omega)]^d:  {\bm v} = {\bm u}_D \mbox{ on }\partial\Omega_D \right\}, \el \\
{\bm V}_0 &:=& \left\{{\bm v} \in [H^1(\Omega)]^d:  {\bm v} = \boldsymbol{0} \mbox{ on }\partial\Omega_D \right\}. \el
\end{eqnarray}
The variational form of \eqref{NS-1}--\eqref{BC-N} is given by: find $({\bm u},p)\in {\bm V} \times L^2(\Omega)$, with ${\bm u}$ satisfying the initial condition \eqref{IC}, such that
\begin{align}
&\int_{\Omega} \frac{\partial{\bm u}}{\partial t}\cdot{\bm v}\de\mathbf{x}+\int_{\Omega}\left({\bm u}\cdot\nabla {\bm u}\right)\cdot{\bm v}\de \mathbf{x}
+\nu \int_{\Omega}\nabla {\bm u}\cdot \nabla{\bm v}\de \mathbf{x}- \int_{\Omega}p\nabla \cdot{\bm v}\de\mathbf{x}\cl
&\hspace{3cm} =\int_{\Omega}{\bf f}\cdot{\bm v}\de\mathbf{x}  + \int_{\partial \Omega_N}{\bf g}\cdot{\bm v}\de\mathbf{x} ,
\qquad\forall\,{\bm v} \in {\bm V}_0, \label{eq:weakNS-1}\\
& \int_{\Omega}q\nabla \cdot{\bm u}\de\mathbf{x} =0, \qquad\forall\, q \in L^2(\Omega).  \label{eq:weakNS-2}
\end{align}

The nonlinearity in problem \eqref{eq:weakNS-1}-\eqref{eq:weakNS-2} can produce a loss of uniqueness of the solution, with multiple solutions branching from a known solution at a bifurcation point. Specifically, Sec.~\ref{sec:cav} considers a cavity model undergoing a Hopf bifurcation, while Sec.~\ref{sec:cha} considers a channel model undergoing a Pitchfork bifurcation.

\subsection{Rayleigh-B\'{e}nard cavity flow}
\label{sec:cav}

\begin{figure}[ht!]
\begin{center}
 \includegraphics[scale=1]{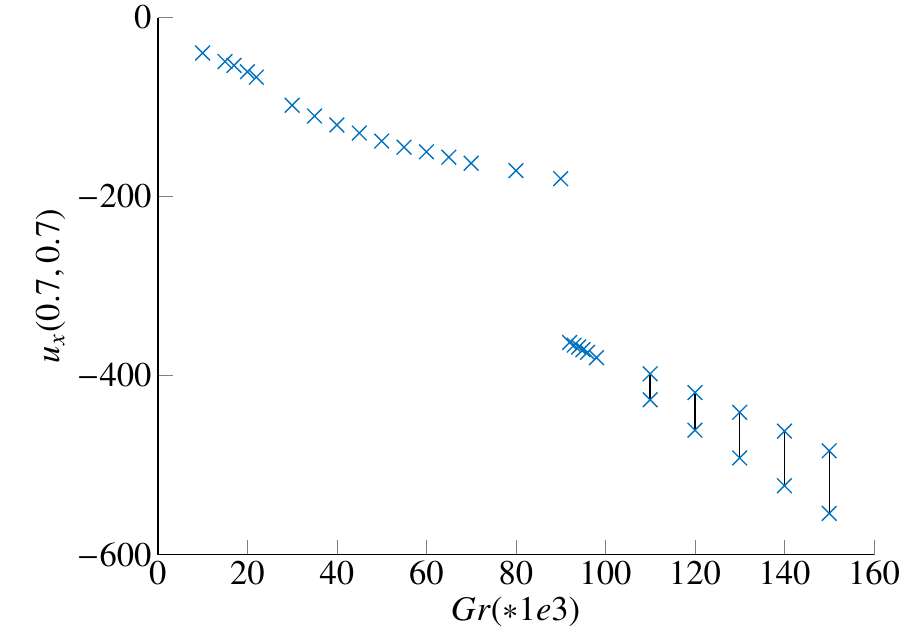} 
 \caption{Bifurcation diagram for the cavity flow. For Gr $> 100\mathrm{e}{3}$, the flow is time-periodic and the range of values is indicated.}
 \label{hess:cavity_bif_dia}
\end{center}
\end{figure}

The Rayleigh-B\'{e}nard cavity flow considered here was 
introduced in~\cite{Roux:GAMM} and has been studied in various works
(see, e.g. \cite{Gelfgat:Ref11,Hess2019CMAME,PR15}).
The computational domain $\Omega$ is a rectangular cavity with height $1$ and length $4$. A no-slip  boundary condition (i.e., a homogeneous Dirichlet condition) is prescribed on the entire boundary and body forcing is given by:
\begin{align}
    {\bm f} = (0, \text{Gr}\nu^2 x)^T, \label{eq_f_cavity}
\end{align}
where $x$ is the horizontal coordinate and Gr is the Grashof number, i.e., the dimensionless number that approximates the ratio of the buoyancy to viscous force acting on a fluid. 
Problem \eqref{NS-1}--\eqref{BC-N} with forcing term \eqref{eq_f_cavity} undergoes several Hopf-type bifurcations as the Grashof number, which is our model parameter, varies in the interval $[10\mathrm{e}{3}, 150\mathrm{e}{3}]$. 
Reference solutions can be found in \cite{Hess2019CMAME,HessQuainiRozza_ACOM_2022} for even larger variations of the Grashof number. 
The bifurcation diagram is shown in Fig.~\ref{hess:cavity_bif_dia}:
it shows the value of the horizontal velocity at the point $(0.7, 0.7)$ under varying Grashof number. 
As is explained in \cite{Hess2019CMAME}, there exists steady single roll solutions until a Grashof number of about Gr=30$\mathrm{e}{3}$, double roll solutions until about Gr=90$\mathrm{e}{3}$ and three roll solutions until about Gr=100$\mathrm{e}{3}$. 
For higher Grashof numbers, we observe time-periodic solutions, such that a whole range of values is indicated.


For the purpose of model reduction, we will work with the trajectory at a fixed parameter, i.e., at a given Grashof number. By trajectory we mean the set of vectors
\begin{equation}
 \{ \mathbf{ x }^{k} \in \mathbb{R}^\mathcal{N} \quad | \quad  k = 1, \ldots, m \},
 \label{hess:define_trajectory}
\end{equation}
that correspond to the discretized velocity ${\bm u}$ at the time steps $\{t_k | k = 1, \ldots, m\}$. A component of $\mathbf{ x }^{k}$ is denoted $ x^{k}_i$ and the trajectory of a single component is $x_{i}$. In \eqref{hess:define_trajectory}, $\mathcal{N}$ denotes the size of the spatial discretization.

In \cite{HessQuainiRozza_ACOM_2022}, we introduced
a multi-stage ROM approach, consisting of proper orthogonal decomposition (POD), dynamic mode decomposition (DMD) and manifold interpolation. 
In this paper, we extend such approach to cover trajectories from a common initial value and we add stabilization techniques to the DMD algorithm.

\subsection{Channel flow}
\label{sec:cha}

\begin{figure}[ht!]
\begin{center}
 \includegraphics[scale=.2]{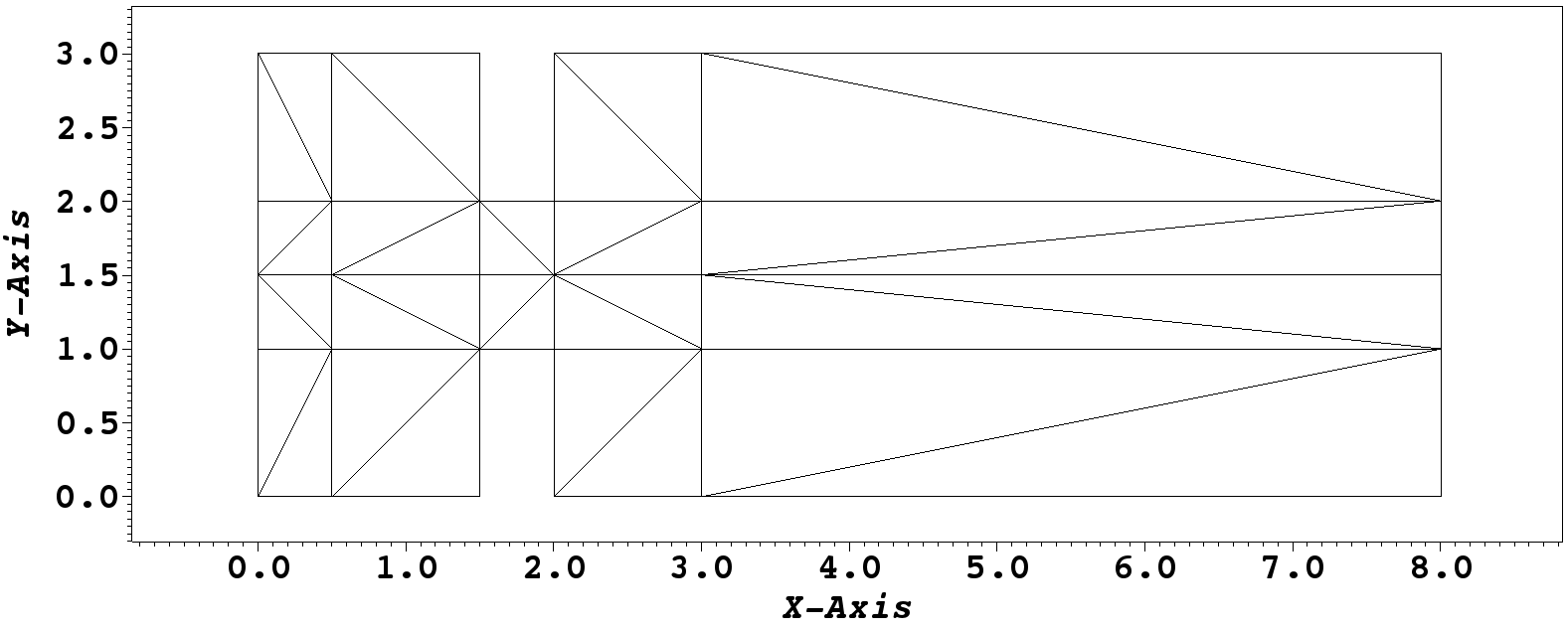} 
 \caption{Computational reference domain for the channel flow at narrowing width $w=1.0$.}
 \label{hess:domain_channel}
\end{center}
\end{figure}

We consider a channel flow through a narrowing of varying width $w$. 
The computational domain is shown in Fig.~\ref{hess:domain_channel} 
and is based on 36 triangular elements. The spectral element solver 
Nektar++\footnote{See \texttt{www.nektar.info}} is used to compute full-order solutions.
This benchmark test has been adapted from various papers, see, e.g., \cite{fearnm1,drikakis1,hawar1,mishraj1} and the references cited therein.
At the inlet, We prescribe a parabolic horizontal velocity component with maximum $9/4$ and 
zero vertical component, as shown in \cite{HessQuainiRozza2022_ETNA}. 
At the outlet, we prescribe a homogeneous Neumann condition (i.e., $\bf g = 0$ in \eqref{BC-N}). On the rest of the boundary, we impose a no-slip condition. 

We vary $\nu$ in the interval $[0.1, 0.2]$ and $w$ in the interval $[0.5, 1.0]$. In those intervals, problem \eqref{NS-1}--\eqref{BC-N} with zero forcing term
undergoes a supercritical pitchfork bifurcation. 

\begin{figure}[ht!]
\begin{center}
 \includegraphics[scale=.9]{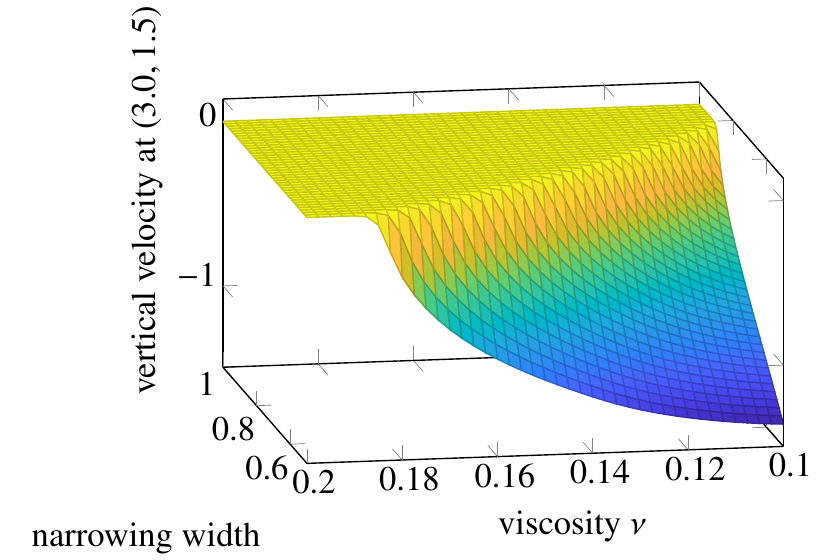} 
 \includegraphics[scale=.9]{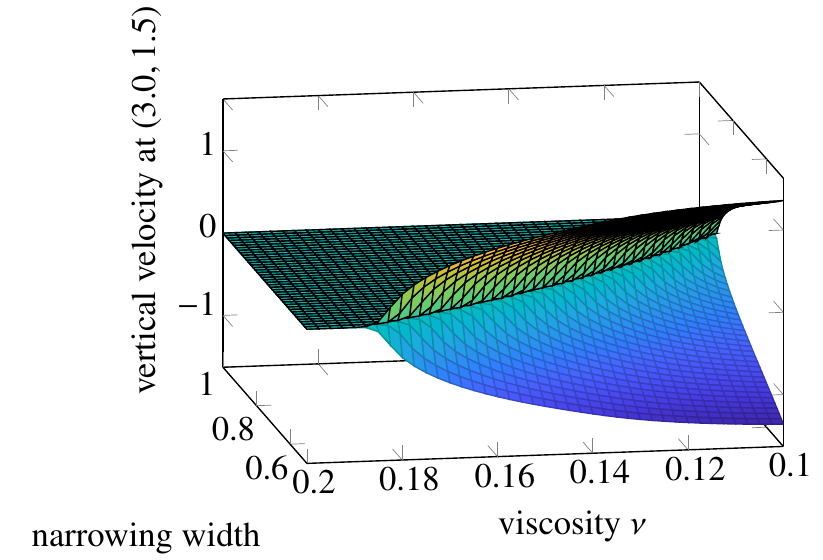} 
 \caption{Left: Bifurcation diagram for the channel flow showing the lower solution branch. Right: Bifurcation diagram showing both solution branches.}
 \label{hess:both_FOM}
\end{center}
\end{figure}

Reference solutions are plotted in \cite{HessQuainiRozza2022_ETNA} on a fine $40 \times 41$ uniform parameter grid. This grid will also serve to assess the accuracy of different methods.
The vertical velocity at point $(3.0, 1.5)$ serves to visualize the bifurcation in a bifurcation diagram, see Fig.~\ref{hess:both_FOM} (left).
Here, the model reduction approach presented in \cite{HessQuainiRozza2022_ETNA} will be extended to capture multiple solution branches as shown in Fig.~\ref{hess:both_FOM} (right).

\section{Reduced order model for the Rayleigh-B\'{e}nard cavity flow}\label{sec:cavity}

We will start this section with a brief recapitulation of the model reduction approach introduced in \cite{HessQuainiRozza_ACOM_2022}, 
which is a multi-step procedure consisting of proper orthogonal decomposition (POD), dynamic mode decomposition (DMD) and manifold interpolation. 
Then, we will present improvements in Sec.~\ref{sec:IS} and \ref{sec:stab}.

The parameter domain, i.e., Gr interval $[100\mathrm{e}{3}, 150\mathrm{e}{3}]$, is sampled uniformly and the full order time trajectories of the flow velocities are computed at each sampled Grashof number. 
A projection space for all trajectories in the parameter domain is computed from the sampled trajectories through POD. 
Let $N$ be the dimension of such projection space.
The POD provides a projected trajectory at each sampled Gr, which is used as input for the dynamic mode decomposition (DMD). 
The DMD computes a linear operator $A \in \mathbb{R}^{N \times N}$, which approximates the dynamics as
\begin{equation}
 \mathbf{ x }^{k+1} \approx A \mathbf{ x }^k \quad \forall k = 1, \ldots, m-1.
 \label{req_DMD}
\end{equation}
Matrix $A$ is called Koopman operator. 
It is decomposed in the DMD modes $U_r \in \mathbb{R}^{N \times r}$ and 
reduced Koopman operator $A_r \in \mathbb{R}^{r \times r}$ with $r \leq N$ as
\begin{equation}
 A = U_r A_r U_r^T .
 \label{Hess:decomposition_Koopman_operator}
\end{equation}
Finally, during the online phase we apply 
manifold interpolation to obtain the reduced DMD operator needed to evaluate the trajectory at a 
new parameter of interest. Once the new reduced Koopman operator is computed with manifold interpolation, we use \eqref{Hess:decomposition_Koopman_operator} to get 
the Koopman operator.

\subsection{Improvement I: Using a common initial state}\label{sec:IS}

The approach presented in \cite{HessQuainiRozza_ACOM_2022} requires an initial state $\mathbf{ x }^1$ that is on the limit cycle or at least close to the limit cycle. This has two disadvantages. First, the initial state is different for each parameter and it requires some previous time advancement of the full order model. Second, it is also not trivial to automate the detection of the limit cycle, which makes manual inspection often necessary. Thus, the model reduction approach could be much improved if a common starting value for all parameters of interest could be used.
 
\begin{figure}[htb!]
\begin{center}
 \includegraphics[width=.6\textwidth]{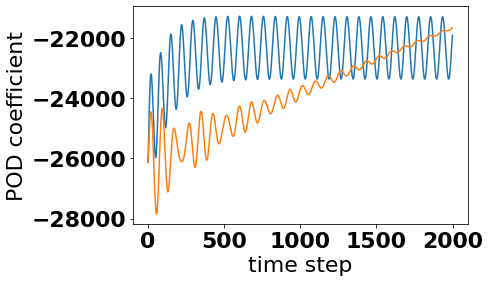} 
\caption{Time evolution of the coefficient of the first POD mode for the full-order model (blue) at Gr$=120\mathrm{e}{3}$ and its standard DMD approximation (orange).  
The initial value for the DMD is on the limit cycle at Gr$=150\mathrm{e}{3}$.
}
 \label{hess:std_DMD_common_init}
\end{center}
\end{figure}

Fig.~\ref{hess:std_DMD_common_init} shows in orange the time evolution of the coefficient of the first POD mode
at Gr$=120\mathrm{e}{3}$ when the initial value is on the limit cycle at Gr$=150\mathrm{e}{3}$. 
Since the provided initial value is far from the limit cycle, the DMD approximation turns out to be highly inaccurate. If
the initial value was given closer the limit cycle at Gr$=120\mathrm{e}{3}$ (for example, after $500$ time steps with respect to the initial value used for Fig.~\ref{hess:std_DMD_common_init}), the DMD would produce an accurate approximation, as is shown in \cite{HessQuainiRozza_ACOM_2022}.

\begin{figure}[htb!]
\begin{center}
 \includegraphics[width=.6\textwidth]{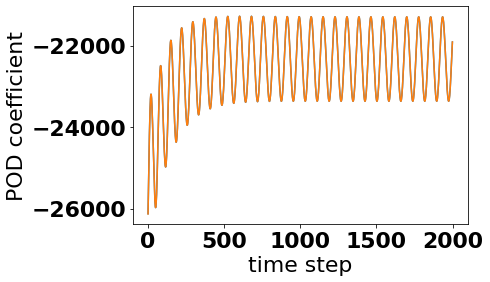} 
\caption{Time evolution of the coefficient of the first POD mode in time of the full-order model (blue) at Gr$=120\mathrm{e}{3}$ and its Hankel-DMD approximation (orange). The initial value for the DMD is on the limit cycle at Gr$=150\mathrm{e}{3}$.
The two curves overlap, which is why only the orange curve can be seen.}
 \label{hess:Hankel_DMD_common_init}
\end{center}
\end{figure}

Hankel-DMD \cite{HankelDMD_Arbabi2017}, one of
several variants of DMD that exist in the literature,
is able to capture the trajectory from different initial values. See Fig.~\ref{hess:Hankel_DMD_common_init}. 
The Hankel-DMD constructs Hankel matrices for each POD coefficent $x_i$.
In \cite{HankelDMD_Arbabi2017}, the coefficients $x_i$ are called observables.

Let 

\begin{eqnarray}
H^i = 
\left[
\begin{array}{cccc}
x^0_i & x^1_i & \ldots & x^n_i  \\
x^1_i & x^2_i & \ldots & x^{n+1}_i \\
\vdots & \vdots & \ddots & \vdots \\
x^{m-1}_i & x^{m}_i & \ldots & x^{m+n-1}_i 
\end{array}
\right]
\label{Hess:Hankel}
\end{eqnarray}
and $A H^i$ denote the Hankel matrix shifted by one time step.
The Hankel matrices are then scaled and arrayed in composite matrices. The exact DMD is  invoked to compute the Koopman eigenmodes and eigenfunctions. See Algorithm 4 in \cite{HankelDMD_Arbabi2017} for details. Then, in principle, it would be possible to also apply manifold interpolation to the matrices of Koopman eigenmodes and the Koopman eigenfunctions.

\subsection{Improvement II: Stabilizing the dynamic mode decomposition}\label{sec:stab}

\begin{figure}[htb!]
\begin{center}
 \includegraphics[width=.6\textwidth]{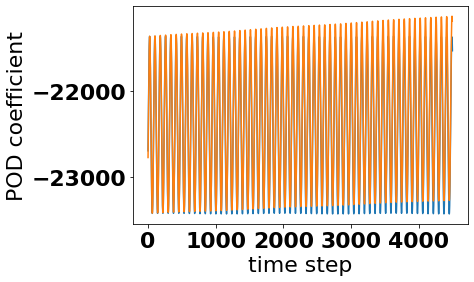} 
 \caption{
Evolution over a long time interval of the coefficient of the first POD mode for the full-order model (blue) at Gr$=120\mathrm{e}{3}$ and its DMD approximation (orange).}
 \label{hess:Gr120_long_time}
\end{center}
\end{figure}

When the DMD approximation is evaluated over a long time interval, the accuracy of the approximation tends to decrease. This is due to the fact that typically the eigenvalues of the Koopman matrix $A$ are not exactly on the complex unit circle. 
See Fig.~\ref{hess:Gr120_long_time}, which reports the DMD approximation of the coefficient of the first POD mode at Gr$=120\mathrm{e}{3}$ from an initial value close to the limit cycle. 
If an eigenvalue of the Koopman matrix has magnitude larger than one, then the corresponding DMD mode will be amplified and the trajectory will approach arbitrary large values for a sufficiently large number of time steps. 
On the other hand, if all eigenvalues have magnitude smaller than one, then all DMD modes will be dampened and the trajectory will approach zero. 

\begin{figure}[htb!]
\begin{center}
 \includegraphics[width=.6\textwidth]{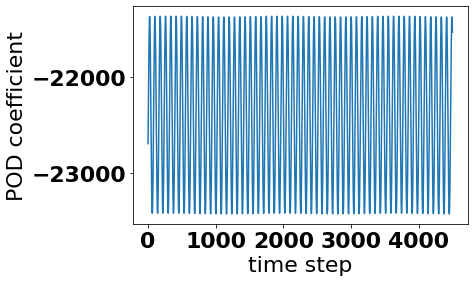} 
 \caption{
 Evolution over a long time interval of the coefficient of the first POD mode at Gr$=120\mathrm{e}{3}$ computed by the stabilized DMD.}
 \label{hess:Gr120_stab_long_time}
\end{center}
\end{figure}

In order to avoid this issue, we perform a simple scaling of the eigenvalues to magnitude one using the eigendecomposition of Koopman matrix $A$ and then reconstructing the stabilised Koopman matrix. The resulting long time trajectory at Gr$=120\mathrm{e}{3}$ is shown in Fig.~\ref{hess:Gr120_stab_long_time}. We observe that it is stable and accurate as it keeps constant amplitude and phase over a long time period of time.

\section{Reduced order model for the channel flow}

In \cite{HessQuainiRozza2022_ETNA}, we proposed a
localized model order reduction method that constructs several local bases, each of which is used for parameters belonging to a different subregion of the parameter domain. We paid particular attention to the criterion used to 
assign a local ROM to each subregion of the parameter domain.
Snapshots are collected over the whole parameter domain and the k-means algorithm is used to cluster the snapshots.
A POD of the clustered snapshots determines a projection space for each cluster, which in turn associates a ROM to each cluster. In order to find a mapping from a parameter location to a cluster, we trained an artificial neural network (ANN). Then, this mapping assigns a cluster to each online parameter of interest. 
The ANN training strategy can make a significant difference for the overall accuracy. In this work, we will show our localized ROM approach can be used also to resolve multiple solutions at a single parameter location.

When one is interested in multiple solution branches, continuation and deflation techniques  \cite{PintorePichiHessRozzaCanuto2021} allow to follow the various branches. In the range of Reynolds number we consider in this paper, only up to two stable solutions are present at any parameter point. There exist also a branch of unstable solutions (\cite{PintorePichiHessRozzaCanuto2021,Hess2019CMAME}), but it is not trivial to compute them as they will converge to one of the stable solutions under small perturbations. Thus, we prefer to focus on the stable solution branches.

\begin{figure}[htb!]
\begin{center}
 \includegraphics[scale=.2]{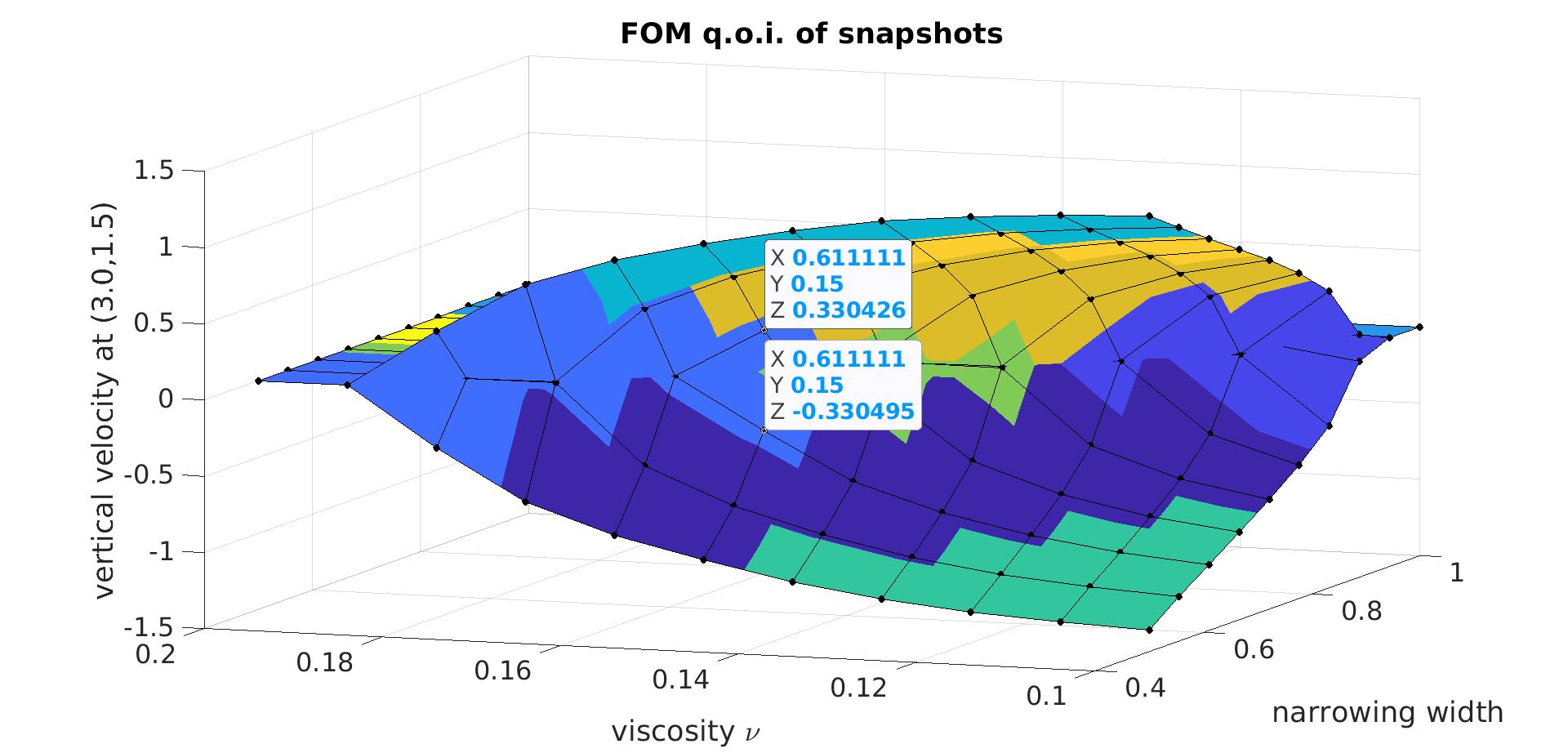}
\caption{Division of the snapshots into ten clusters
according to the k-means algorithm. A different color is used for each cluster. The light blue cluster includes 
snapshots from both the upper and lower branches at the same parameter point, as highlighted also by the markers.}
 \label{hess:cluster_upper_lower_naive}
\end{center}
\end{figure}

We consider a uniform $10 \times 11$ sampling of the parameter domain, with one or two solutions being computed depending on whether at a parameter point the solution is unique or not.
The collected snapshots are clustered by the k-means algorithm, with the appropriate number of clusters determined from the k-means energy. Alternatively, the silhouette score  can be used to determine the number of clusters, see \cite{Hess2019CMAME}.
Fig.~\ref{hess:cluster_upper_lower_naive} shows the division  of the snapshots into ten clusters, each represented with a different color. For each cluster, the POD computes projection spaces in the same way as when only single solutions exist for a given parameter.

\begin{figure}[hbt!]
\begin{center}
 \includegraphics[scale=1]{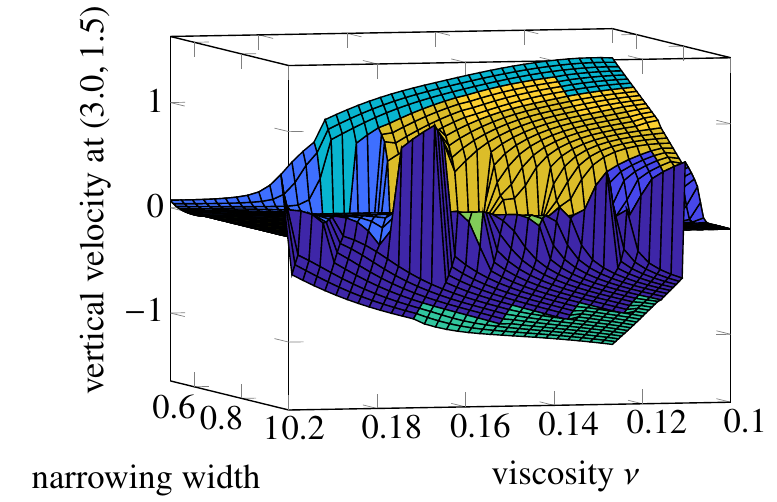} 
\caption{Reconstruction of the bifurcation diagram 
with the local ROM approach that uses the
distance to snapshot location criterion.}
 \label{hess:resolve_upper_lower_naive}
\end{center}
\end{figure}

The first criterion we consider to select the local ROM for an online parameter of interest is the closest snapshot parameter location. This criterion entails finding the snapshot parameter location closest to the online parameter so that one uses the projection space corresponding to the cluster the snapshot belongs to. Fig.~\ref{hess:resolve_upper_lower_naive} shows the reconstruction of the bifurcation diagram according to this selection criterion. We observe that close to the onset of the bifurcation the diagram is poorly reconstructed
as only a single solution is captured at parameter locations where two solutions exists. This issue could have been anticipated from looking at the clusters in Fig.~\ref{hess:cluster_upper_lower_naive} because 
some clusters, e.g., the light blue one (with the two markers), include snapshots from both the upper and lower branches. Then, when solving within the associated projection space only one solution is found. We obtain a mean relative $L^2$ error for the velocity of $5\%$.
By computing the optimal cluster selection, we find that with this clustering the mean error is bounded from below by $0.41\%$. 

\begin{figure}[htb!]
\begin{center}
 \includegraphics[scale=.2]{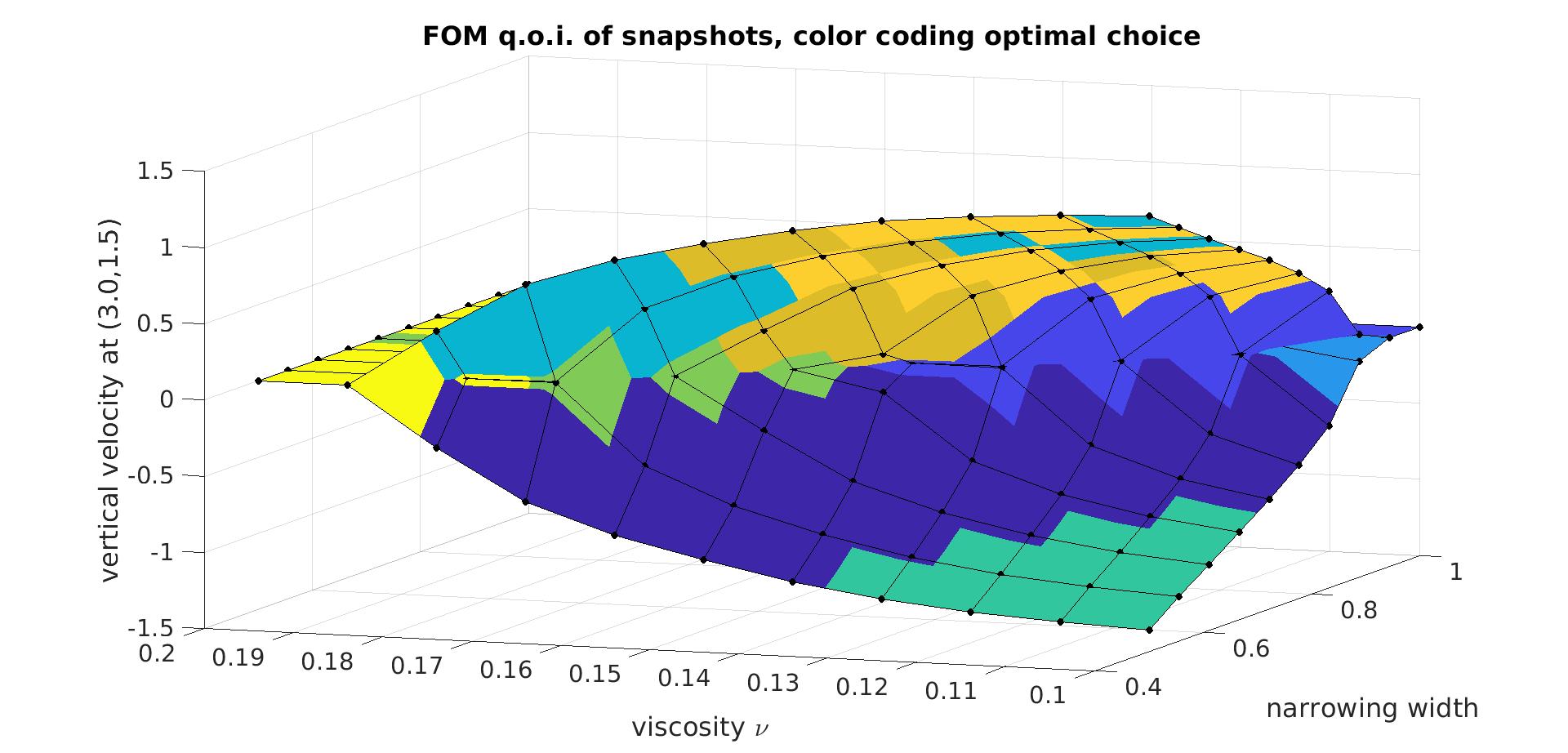} 
 \caption{
Selection of the cluster that gives the lowest error at each snapshot location.
The colors refer to the clusters in Fig.~\ref{hess:cluster_upper_lower_naive}.}
 \label{hess:cluster_upper_lower_smart}
\end{center}
\end{figure}

The poor resolution of the bifurcation diagram in Fig.~\ref{hess:resolve_upper_lower_naive} could be avoided by employing a deflation scheme within each local ROM \cite{PintorePichiHessRozzaCanuto2021}. 
However, another strategy is possible.
After computing the clustering, one solves each local ROM at each snapshot location and collects the relative errors. 
Since this is part of the offline phase, the additional computational effort is negligible.
Fig.~\ref{hess:cluster_upper_lower_smart} shows the resulting cluster selection, which gives the lowest error at each snapshot location. 
A crucial difference with this cluster selection criterion with respect to the parameter clustering in Fig.~\ref{hess:cluster_upper_lower_naive} is that the upper and the lower branches are associated to different clusters.
We train an ANN with the mapping from parameter location  of a snapshot to the local ROM with the lowest error. 
Since some parameters are associated with two snapshots, two such ANNs are computed. 
Each ANN is modeled as a multi-layer perceptron with 5 layers of 512 nodes each and a ReLU activation function. If the online parameter of interest belongs to the part of the parameter domain where a single solution exist, then both ANNs will choose the same cluster and thus local ROM. 
Otherwise, the two ANNs will choose different clusters since, as explained above, no cluster extends to both bifurcation branches.

\begin{figure}[htb!]
\begin{center}
 \includegraphics[scale=1]{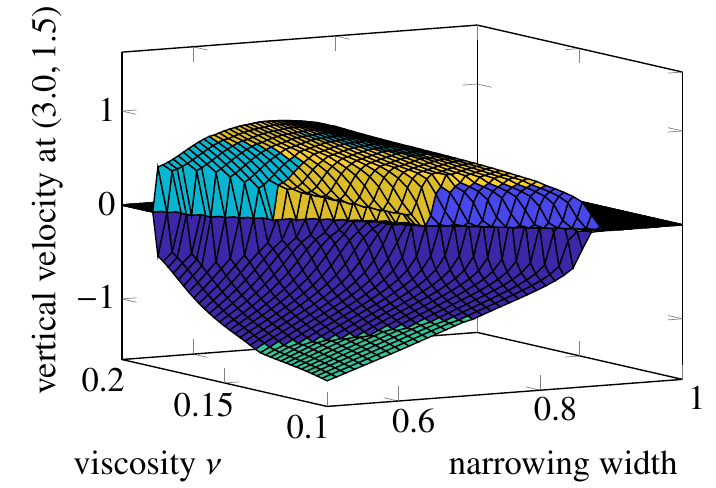} 
 \caption{Reconstruction of the bifurcation diagram 
with the local ROM approach that uses the the ANN cluster selection criterion. The training data are the
relative errors of the local ROMs.}
 \label{hess:resolve_upper_lower_smart}
\end{center}
\end{figure}

The bifurcation diagram obtained with this alternative procedure is shown in Fig.~\ref{hess:resolve_upper_lower_smart}. 
The mean relative $L^2$ error for the velocity is $0.48\%$, which is very close to the $0.41\%$ given by the optimal cluster selection. 
Although the bifurcation diagram in Fig.~\ref{hess:resolve_upper_lower_smart} is an accurate reconstruction of the diagram computed with the FOM reported in Fig.~\ref{hess:both_FOM}, it is not as smooth. 
This is due to the fact that the ANN cluster selection training was performed with respect to the relative $L^2$ error instead of the error at the single degree of freedom shown in the bifurcation diagram. 

\section{Summary and Outlook}

In this paper, we showed how data-driven techniques can aid projection methods to reconstruct bifurcation diagrams
accurately.
Two distinct contexts are considered: cavity flow with a Hopf bifurcation as the Grashof number is varied and jet flow with a supercritical pitchfork bifurcation as the Reynolds number and the jet inlet width change.
For the cavity flow, we proposed to use the Hankel-DMD
algorithm, which allows to approximate the solution at different Grashof numbers starting from the same initial state, unlike the standard DMD algorithm, which requires an ad-hoc initial value. 
In addition, we stabilize the DMD in order to correctly recover the projected cavity flow over long time trajectories.
For the jet flow, we proposed to select the local 
ROM for a new parameter point with an ANN trained using the known relative errors at the snapshot locations. With this, our local
ROM approach manages to recover multiple solutions. 

Future work could go towards combining the Hankel-DMD with manifold interpolation and the local ROM approach towards more complex applications.

\section*{Acknowledgments}

We acknowledge the support provided by the European Research Council Executive Agency by the Consolidator Grant project AROMA-CFD ``Advanced Reduced Order Methods with Applications in Computational Fluid Dynamics" - GA 681447, H2020-ERC CoG 2015 AROMA-CFD, PI G.~Rozza, and INdAM-GNCS 2019-2020 projects.
This work was also partially supported by US National Science Foundation through grant DMS-1953535. A.~Quaini acknowledges 
support from the Radcliffe Institute for Advanced Study at Harvard University where she has been the 2021-2022 William and Flora Hewlett Foundation Fellow.

 \bibliographystyle{siam}
 \bibliography{rbsissa, latexbi}

\end{document}